\documentstyle[10pt]{article}

\newcommand{\rom}{\mathrm}

\def\N{{\rom I\kern-.1567em N}}
\def\R{{\rom I\kern-.1567em R}}
\def\C{{\rom C\kern-6.5pt  
          \vrule height 7.7pt width 0.4pt depth -0.5pt \phantom {.}\;}}             

\newtheorem{satz}{Theorem}
\newtheorem{lemma}[satz]{Lemma}
\newtheorem{proposition}[satz]{Proposition}
\newtheorem{coro}[satz]{Corollary}
\newcommand{\eps}{\varepsilon}
\newcommand{\falle}{\;\;\forall}

\newcommand{\gen}{\rightarrow}

\newcommand{\nach}[1]{\stackrel{#1}\rightarrow}
\newcommand{\aus}{\subset}
\newcommand{\Norm}[1]{\Bigl\|#1\Bigr\|}
\newcommand{\norm}[1]{\|#1\|}
\newcommand{\Betr}[1]{\Bigl| #1  \Bigr|}
\newcommand{\betr}[1]{| #1  |}
\newcommand{\eing}[1]{|_{#1}}
\newcommand{\ebew}{\hfill{\rule{1.2ex}{1.2ex}}}
\newcommand{\bgl}{\begin{eqnarray}}
\newcommand{\bglst}{\begin{eqnarray*}}
\newcommand{\egl}{\end{eqnarray}}
\newcommand{\eglst}{\end{eqnarray*}}

\newcommand{\folgt}{\Rightarrow}
\newcommand{\leins}{l^{1}}
\newcommand{\lunendl}{l^\infty}
\newcommand{\Ref}[1]{(\ref{#1})}
\newcommand{\gdw}{\Leftrightarrow}
\newcommand{\eins}{{1}}

\newcommand{\charF}[1]{\eins_{#1}}

\newcommand{\gentau}{\nach{\tau}}

\newcommand{\mdE}{|\;}
\newcommand{\LL}{{\rom L}}
\newcommand{\Cstalg}{C$^{*}$-al\-ge\-bra}
\newcommand{\Lunendl}{{\rom L}^{\infty}}
\newcommand{\Leins}{{\rom L}^{1}}
\newcommand{\Lnull}{{\rom L}^{0}}
\newcommand{\Ree}{{\rom R}{\rom e}\,}

\begin{document}
\begin{center}
\begin{bf}Perturbation of $\leins$-copies and measure convergence
in preduals of von Neumann algebras\end{bf}
\medskip\\ H. Pfitzner\end{center}
\begin{abstract}\noindent
Let $\Leins$ be the predual of a von Neumann algebra with
a finite faithful normal trace.
We show that a bounded
sequence in $\Leins$ converges to $0$ in measure if and only
if each of its subsequences admits another subsequence which
converges to $0$ in norm or spans $\leins$ "almost isometrically".
Furthermore we give a quantitative version of an
essentially known result concerning the perturbation of a sequence
spanning $\leins$ isomorphically in the dual of a C$^*$-algebra.
\end{abstract}
\bigskip\bigskip\bigskip
{\bf\S 1 Introduction, main results}
\bigskip\\
The present article deals with convergence in probability
in $\Leins$-spaces from a functional analytic point of view.
The $\Leins$-spaces in question are the preduals of
von Neumann algebras with finite faithful normal traces.
To consider an easy example we look at the commutative case:
Let $(\Omega,\Sigma,\mu)$ be a finite measure space,
let $(f_n)$ be a bounded sequence in $\Leins(\Omega,\Sigma,\mu)$.
If (appropriately chosen representatives of)
the $f_n$ have pairwise disjoint supports then clearly
$(f_n)$ converges to $0$ in measure.
From the functional analytic point of view such a sequence,
up to normalization, is the canonical
basis of an isometric copy of $\leins$.
If one perturbes $(f_n)$ by a norm null sequence $(g_n)$ then $(f_n+g_n)$
still $\mu$-converges to $0$ and spans $\leins$ almost isometrically
(in a sense to be made precise below in \S 2).
It has been known \cite[Th.\ 2]{KadPel}
(see also \cite[Th.\ 3, Rem.\ 6bis]{Pi-bases}) for quite a time
that, roughly speaking,  these are essentially the only examples
of $\mu$-null sequences.

Theorem \ref{theo_L1}
contains the analogous statement for the predual of a
von Neumann algebra with finite faithful normal trace.
(For notation and definitions see \S 2.)
\begin{satz}\label{theo_L1}
Let $(x_n)$ be a bounded sequence in $\Leins ({\cal N},\tau)={\cal N}_*$
where $({\cal N},\tau)$ is a
von Neumann algebra with a finite normal faithful trace $\tau$.
Then the following  assertions are equivalent.
\begin{description}
\item{\makebox[1.8em]{(i)}}
$x_n\gentau 0$.
\item{\makebox[1.8em]{(ii)}}
For each subsequence $(x_{n_{k}})$ of $(x_n)$
there are a subsequence $(x_{n_{k_l}})$ and
a sequence $(y_l)$ of pairwise orthogonal elements
of $\Leins ({\cal N},\tau)$ such that $\norm{x_{n_{k_l}}-y_l}_1\gen 0 $.
\item{\makebox[1.8em]{(iii)}}
For each subsequence $(x_{n_k})$ of $(x_n)$
there is a subsequence $(x_{n_{k_l}})$
which tends to $0$ in $\norm{\cdot}_1$ or 
spans $\leins$ almost isometrically.
\item{\makebox[1.8em]{(iv)}}
For each subsequence $(x_{n_k})$ of $(x_n)$
there is a subsequence $(x_{n_{k_l}})$
which tends to $0$ in $\norm{\cdot}_1$ or 
spans $\leins$ asymptotically.
\end{description}
The implications (i) $\gdw$ (ii) $\folgt$ (iii) $\gdw$ (iv) hold also
for unbounded sequences $(x_n)$,
the implications (iii) $\folgt$ (ii), (i) do not.
\end{satz}

Implication (i)$\folgt$(ii) has already appeard as
a special case of a result of Sukochev \cite[Prop.\ 2.2]{Suk}.
The other nontrivial implication
(iii)$\folgt$(ii) follows immediately from Theorem
\ref{prop_leins_asy} which holds for the predual
of any von Neumann algebra and is of independent interest:
\begin{satz}\label{prop_leins_asy}
Let ${\cal N}$ be an arbitrary von Neumann algebra and $(\phi_m)$ a bounded
sequence in its predual ${\cal N}_*$.
If $(\phi_m)$ spans $\leins$ almost isometrically
then there are a subsequence $(\phi_{m_l})$ of
$(\phi_m)$ and
a sequence $(\tilde{\phi}_l)$ of pairwise orthogonal functionals
in ${\cal N}_*$ such that
$\norm{\phi_{m_l} -\tilde{\phi}_{l}} \gen 0$ as
$l\gen\infty$.\\
This amounts to saying that
there are pairwise orthogonal projections $s_l$ and
pairwise orthogonal projections $t_l$ in ${\cal N}$ such that
$\norm{\phi_{m_l}-t_l \phi_{m_l} s_l}\gen 0$ as
$l\gen\infty$.
\end{satz}
It is natural to ask what can be improved in Theorem
\ref{prop_leins_asy} if one replaces the predual of the von Neumann algebra
by the dual of a \Cstalg. At the time of this writing this is not
clear. What we have among other things is
\begin{proposition}\label{theo_A1}
Let $(\phi_m)$ be a bounded sequence that spans $\leins$
almost isometrically in the dual of an arbitrary
\Cstalg\ $A$.
Then, given $\eps >0$, there are pairwise orthogonal positive normalized
elements $(a_n)$ and pairwise orthogonal positive normalized
elements $(b_n)$ in $A$ such that
$\norm{\phi_{m_n}-b_n \phi_{m_n} a_n}<\eps$ for an appropriate
subsequence $\phi_{m_n}$ and all $n\in\N$ .
\end{proposition}
For a more detailed discussion see \S 6.\bigskip\\
As to the organization of the paper, after recalling
some notation and definitions in the
next section we gather some auxiliary results in \S 3
in order to prove Theorem \ref{prop_leins_asy} in \S 4.
In \S 5 we  prove Theorem \ref{theo_L1} for the sake of completeness
although, as already mentionned, it follows essentially from
\cite[Prop.\ 2.2]{Suk} and Theorem \ref{prop_leins_asy}.
In \S 6 perturbations of $\leins$-copies in the dual
of C$^*$-algebras are considered.
\bigskip\bigskip\\
{\bf\S 2 Notation, definitions}\\
Let $(x_n)$ be a sequence of nonzero elements in a Banach space $X$.

We say that
\begin{em}$(x_n)$ spans $\leins$ $r$-isomorphically\end{em}
or just \begin{em}isomorphically\end{em}
if there exists $r>0$ (trivially $r\leq 1$) such that
$r(\sum_{n=1}^{\infty}\betr{\alpha_n})\leq
\norm{\sum_{n=1}^{\infty} \alpha_n \frac{x_n}{\norm{x_n}}}   \leq
        \sum_{n=1}^{\infty}\betr{\alpha_n}$
for all scalars $\alpha_n$ (the second inequality being trivial).

We say that
\begin{em}$(x_n)$ spans  $\leins$ almost isometrically\end{em}
if there is a sequence $(\delta_m)$ in $[0,1[$
tending to $0$ such that
$(1-\delta_m)\sum_{n=m}^{\infty}\betr{\alpha_n}\leq
\norm{\sum_{n=m}^{\infty} \alpha_n \frac{x_n}{\norm{x_n}}}   \leq
         \sum_{n=m}^{\infty}\betr{\alpha_n}$
for all $m\in\N$.\smallskip

Trivially the property of spanning $\leins$ almost isometrically
passes to subsequences.
Recall that James' distortion theorem
(see \cite{Jam64} or \cite{Die-Seq})
for $\leins$ says that every isomorphic
copy of $\leins$ contains an almost isometric copy
of $\leins$. To be more precise,
let $r>0$, $[0,1[\ni\delta_n\gen0$, 
and let  $(x_n)$ be a normalized basis spanning $\leins$
$r$-isomorphically.
Then it follows from the proof of \cite{Jam64} that
there is a sequence $(\lambda_{i})$ of scalars
and a sequence $(F_{n})$ of pairwise disjoint finite
subsets of $\N$ such that
$(1-\delta_m)\sum_{n=m}^{\infty}\betr{\alpha_n}\leq
\norm{\sum_{n=m}^{\infty} \alpha_n {y_n}}   \leq
         \sum_{n=m}^{\infty}\betr{\alpha_n}$
for all scalars $\alpha_n$ and all $m\in\N$
where $y_{n}=\sum_{i\in F_{n}}\lambda_{i}x_i$
and where $\sum_{i\in F_{n}}\betr{\lambda_{i}}\leq \frac{1}{r}$ for
all $n\in\N$.

Finally $(x_n)$ is said {\em to span $\leins$
asymptotically isometrically} or just
{\em to span $\leins$ asymptotically} if there is a sequence
$(\delta_n)$ in $[0,1[$ tending to $0$ such that
$\sum_{n=1}^{\infty}(1-\delta_n)\betr{\alpha_n}\leq
\norm{\sum_{n=1}^{\infty} \alpha_n \frac{x_n}{\norm{x_n}}}   \leq
         \sum_{n=1}^{\infty}\betr{\alpha_n}$
for all scalars $\alpha_n$.
We say that a Banach space is isomorphic
(respectively almost isometric, respectively asymptotically isometric) to
$\leins$ if it has a basis with the corresponding property.
Clearly a sequence spanning $\leins$ asymptotically spans
$\leins$ almost isometrically.
The main result of \cite{DJLT} states that the converse does not
hold because there are almost isometric copies of $\leins$
which do not contain $\leins$ asymptotically.
However, it follows from \cite{Pfi-Fixp} that this cannot happen in
the predual of a von Neumann algebra because each sequence spanning
$\leins$ almost isometrically in a von Neumann predual contains a
sequence spanning $\leins$ asymptotically (cf. (iii) $\folgt$ (iv)
in the proof of Theorem \ref{theo_L1}).
Note that the present definitions of almost and asymptotically isometric
differ slightly from those in \cite{DJLT}, \cite{Pfi-Fixp} by the term
$x_n /\norm{x_n}$ but that, of course, for normalized sequences
the definitions are the same.
Note also the technical detail that because of this term
one might have $\norm{x_n}\gen 0$ for a sequence spanning $\leins$
isomorphically (or almost or asymptotically isometrically)
whereas sequences that are
equivalent to the canonical $\leins$-basis (\cite[p.\ 43]{Die-Seq})
are uniformly bounded away from $0$.\medskip

The dual of a Banach space $X$
is denoted by $X'$.
We work with complex scalars.
Two elements $a,b$ of a  C$^*$-algebra are called orthogonal
- $a\bot b$ in symbols - if $ab^*=0=a^*b$.
Two elements $\phi,\psi$ of the predual of a von Neumann algebra
are called orthogonal - $\phi\bot\psi$ in symbols - if
they have orthogonal right and orthogonal left support projections.
It is well know that $\phi\bot\psi$ if and only if the linear span
of $\phi$ and $\psi$ is isometrically isomorphic to the two-dimensional
$\leins_2$; if $\phi$ and $\psi$ are positive
they are orthogonal if and only if
$\norm{\phi-\psi}=\norm{\phi}+\norm{\psi}$.

Let ${\cal N}$ be a von Neumann algebra, $a\in{\cal N}$,
$\phi\in{\cal N}_*$ then $a\phi$ denotes the normal functional 
${\cal N}\ni x\mapsto\phi(xa)$ and $\phi a$ denotes the normal functional
${\cal N}\ni x\mapsto\phi(ax)$.

Let $\tau$ be a finite faithful normal trace
on a von Neumann algebra ${\cal N}$.
The set $I=\{x\in {\cal N}\mdE \tau(\betr{x})<\infty\}$
is an ideal in ${\cal N}$, can be normed by
$x\mapsto \tau(\betr{x})=:\norm{x}_1$ and its Banach space completion 
is denoted by $\Leins=\Leins({\cal N},\tau)$. It is well-known that $\Leins$
is isometrically isomorphic to the predual ${\cal N}_*$ via the map
$\Leins\ni x\mapsto \phi_x \in {\cal N}_*$
where $\phi_x (y)=\tau(xy)$ for
$y\in {\cal N}$ and where $\tau$ is understood as the (well-defined) extension
of $\tau$ from $I$ to $\Leins$ \cite[V.2.18]{Tak}.
In particular, the multiplication on ${\cal N}\times I$ can be extended
to ${\cal N}\times \Leins$, the map $x\mapsto \phi_x$ respects orthogonality
and one has $\betr{\tau(xy)}\leq \norm{x}_1 \norm{y}_{\infty}$ for 
$x\in\Leins$, $y\in\Lunendl=\Lunendl({\cal N},\tau):={\cal N}$.
More generally one can define $\LL^p({\cal N},\tau)$-spaces, $1\leq p<\infty$,
as the sets of those $x\in\Lnull$ for which
$\norm{x}_p:=\tau(\betr{x}^p)^{1/p}<\infty$ where
$\Lnull=\Lnull({\cal N},\tau)$ is the space of
$\tau$-measurable densely defined (in general unbounded)
operators affiliated with ${\cal N}$ and where $\tau$ is understood
as the extension of $\tau$ from ${\cal N}$ to $\Lnull$.
On $\Lnull$ one defines the measure topology as the translation invariant
topology in which the sets
$\{x\in\Lnull\mdE \exists p\in {\cal N}_{ {proj}}:\, xp\in {\cal N},\, \norm{xp}_{\infty}
\leq\eps,\, \tau(p^{\bot})\leq\delta\}$, $\eps,\delta>0$, form a base of the
zero neighborhoods.
(${\cal N}_{ {proj}}$ denotes the set of projections of ${\cal N}$.)
In this topology, $\Lnull$ becomes a (well-defined) metrizable
complete Hausdorff topological
vector $^*$-algebra and all $\LL^p$ embed injectively in $\Lnull$.
In particular, sum and product are well-defined in $\Lnull$.
All this (for the more general case of a semifinite trace)
can be found for example in \cite{Nel-nonc}, \cite[Ch.\ 1]{Terp}
or \cite{Yea-Lp}.

If a sequence $(x_n)$ in $\Lnull$ converges to $x\in\Lnull$ this is
denoted by $x_n \gentau x$.
In this context Chebyshev's inequality reads
$\tau(\chi_{]\eps,\infty[}(\betr{x}))
\leq \tau(\frac{1}{\eps}\betr{x})=\frac{1}{\eps}\norm{x}_1$ for
$x\in\Leins$ - which means in particular
that the norm topology is finer than
the measure topology induced by $\tau$ -
and from \cite[A48]{Jajte} we know that in accordance with the
commutative case, $x_n\gentau0$ if and only if
$\tau(\chi_{]\eps,\infty[}(\betr{x_n}))\gen0$ as $n\gen\infty$
for all $\eps>0$.

Basic properties and definitions which are not explained here
can be found in \cite{Die-Seq} or in \cite{LiTz1}-\cite{LiTz2}
for Banach spaces and in \cite{Ped}, \cite{Tak} for
C$^*$-algebras.
\bigskip\bigskip\bigskip\\
{\bf\S 3 Some auxiliary results}
\bigskip\\
Let us first state an easy lemma which says that almost
isometric and asymptotically isometric $\leins$-copies
are stable with respect to perturbations by norm null sequences.
\begin{lemma}\label{lem_asy_stabil}
Let $(x_n)$, $(y_n)$ be two sequences in a Banach space $X$ such
that $\inf\norm{x_n}>0$,
$\norm{y_n}\gen 0$ and $x_n +y_n\neq 0$.\\
If $(x_n)$ spans $\leins$ almost isometrically
then so does $(x_n + y_n)$.\\
If $(x_n)$ spans $\leins$ asymptotically
then so does $(x_n + y_n)$.
\end{lemma}
{\em Proof}:
Suppose that $(x_n)$ spans $\leins$ almost isometrically.
For all scalar sequences $(\alpha_n)$ one has
\bglst
\lefteqn{\Norm{\sum_{n=m}^{\infty}\alpha_n 
                             \frac{x_n +y_n}{\norm{x_n +y_n}}}}\\
&\geq&
\Norm{\sum_{n=m}^{\infty}\alpha_n \frac{x_n}{\norm{x_n}}}
     -
\Norm{\sum_{n=m}^{\infty}\alpha_n\Bigl(1-\frac{\norm{x_n}}{\norm{x_n +y_n}}\Bigr) \frac{x_n}{\norm{x_n}}}
    -
\Norm{\sum_{n=m}^{\infty}\alpha_n \frac{y_n}{\norm{x_n +y_n}}}  \\
&\geq&
\Bigl((1-\delta_m)\sum_{n=m}^{\infty} \betr{\alpha_n}\Bigr)
     -
 \Bigl(\sup_{n\geq m} \Betr{1-\frac{\norm{x_n}}{\norm{x_n +y_n}}}
\sum_{n=m}^{\infty} \betr{\alpha_n}\Bigr)
           -
\Bigl(\sup_{n\geq m}\frac{\norm{y_n}}{\norm{x_n +y_n}}
                        \sum_{n=m}^{\infty} \betr{\alpha_n}\Bigr)\\
&=&
(1-\delta_m')\sum_{n=m}^{\infty} \betr{\alpha_n}
\eglst
where
$\delta_m'=\delta_m
+\sup_{n\geq m} \betr{1-\frac{\norm{x_n}}{\norm{x_n +y_n}}}
+\sup_{n\geq m}\frac{\norm{y_n}}{\norm{x_n +y_n}}
\gen 0$ as $m\gen\infty$.
Hence $(x_n +y_n)$ spans $\leins$ almost isomorphically.
The asymptotic case is proved similarly.
\ebew\medskip\\
Lemmas \ref{lem_A4} - \ref{lem_A3bis} seem to be known
and are proved mainly for lack of suitable reference.
(In part they overlap with \cite[Lem.\ 3-5]{Pf-V}.)
\begin{lemma}\label{lem_A4}
Let $A$ be a \Cstalg, $\omega$ a positive functional on $A$
and $a, b$ elements of the unit ball of $A$. Then
\bgl
\norm{a\omega-\omega} &\leq& (2\norm{\omega})^{1/2}\,
                         \betr{\,\,\norm{\omega}-\omega(a)}^{1/2}
                         \label{glA4_1}\\
\norm{\omega a-\omega} &\leq& (2\norm{\omega})^{1/2}\,
                         \betr{\,\,\norm{\omega}-\omega(a)}^{1/2}
\label{glA4_2}\\
\norm{b\omega a-\omega} &\leq& (2\norm{\omega})^{1/2}\,
                         (\betr{\,\,\norm{\omega}-\omega(a)}^{1/2}
                          +\betr{\,\,\norm{\omega}-\omega(b)}^{1/2})
\label{glA4_3}
\egl
\end{lemma}
{\em Proof:}
Let $x\in A$ and $\norm{x}\leq1$.
Set $\gamma=\norm{\omega}-\omega(a)$,
thus $\omega(a^*)=\norm{\omega}-\overline{\gamma}$.
Without loss of generality we assume $\norm{\omega}=1$.
The inequality of Cauchy-Schwarz yields
\bglst
\betr{\omega(x)-a\omega(x)}^2
      &=& 
     \betr{\omega(x(\eins -a))}^2
     \leq \omega(xx^*)\,\omega((\eins -a)^*(\eins-a))\\
      &\leq&
    \omega((\eins -a)^*(\eins-a))
    =\omega(\eins-a) -\omega(a^*-a^*a)\\
     &=& \gamma-(1-\overline{\gamma})+\omega(a^*a)
    \leq
    2\Ree\gamma \leq2\betr{1-\omega(a)}
\eglst
whence \Ref{glA4_1}; \Ref{glA4_2} follows analogously;
\Ref{glA4_3} follows from \Ref{glA4_1}, \Ref{glA4_2} and from
$\norm{\omega-b\omega a}\leq\norm{\omega-b\omega}+\norm{b(\omega-\omega a)}
\leq \norm{\omega-b\omega}+\norm{\omega-\omega a}$
\ebew

\begin{lemma}\label{lem_A3}
Let $A$ be a \Cstalg, $\phi$ a functional on $A$
and $a,b$ in the unit ball of $A$.
Then
\bgl
\norm{\phi-a\betr{\phi}\;}
   &\leq&
({2\norm{\phi}})^{1/2}\,\,\betr{\,\norm{\phi} - \phi^*(a)}^{1/2}   \label{glA3_1}\\
\norm{\betr{\phi}-a\phi}
   &\leq&
({2\norm{\phi}})^{1/2}\,\,\betr{\,\norm{\phi} - \phi(a)}^{1/2}   \label{glA3_2}\\
\norm{b\phi a-\phi}
   &\leq&
({2\norm{\phi}})^{1/2}\,\,
\Bigl(\betr{\,\norm{\phi} - \betr{\phi}(a)} 
     + \betr{\,\norm{\phi} - \betr{\phi^*}(b)}\Bigr)^{1/2}.   \label{glA3_3}
\egl
\end{lemma}
{\em Proof}:
Let $\phi=u\betr{\phi}$ be the polar decomposition of $\phi$.
Then the polar decomposition of $\phi^*$ is
$\phi^*=u^*\betr{\phi^*}$  (cf. the proof of \cite[III.4.2]{Tak}),
we have $\phi=\betr{\phi^*}u$,
$\betr{\phi}=u^*\phi=\betr{\phi}^*=\phi^*u$.
Without loss of generality we assume $\norm{\phi}=1$.\\
Inequality \Ref{glA3_2} follows from 
\bglst
\norm{a\phi-\betr{\phi}\,\,}=\norm{au \betr{\phi}-\betr{\phi}\,}
\stackrel{\Ref{glA4_1}}{\leq}
\betr{2(1-\betr{\phi} (au))}^{1/2}
=
\betr{2(1-\phi(a))}^{1/2}.
\eglst
Replacing $\phi$ by $\phi^*$ we have 
$\norm{a\phi^*-\betr{\phi^*}\;}\leq\betr{2(1-\phi^*(a))}^{1/2}$
whence \Ref{glA3_1} by\\ 
$\norm{a\betr{\phi}-\phi}=\norm{(a\phi^*  -\betr{\phi^*}) u }
\leq
\norm{a\phi^* - \betr{\phi^*}\;}$.
\Ref{glA3_3} follows from
\bglst
\norm{\phi-b\phi a}
&\leq&
\norm{\phi - b\phi} + \norm{b\phi - b\phi a}\\
&=&
\norm{\; \betr{\phi^*} u -b \betr{\phi^*} u}+\norm{bu \betr{\phi}-bu \betr{\phi} a} \\
&\leq&
\norm{\,\betr{\phi^*} - b \betr{\phi^*}\,} 
         + \norm{\,\,\betr{\phi} - \betr{\phi} a}\\
&\stackrel{\Ref{glA4_1}\Ref{glA4_2}}{\leq}&
(2\norm{\phi})^{1/2}\,
(\betr{\,\norm{\phi} - \betr{\phi^*}(b)} 
     + \betr{\,\norm{\phi} - \betr{\phi}(a)})^{1/2}.
\eglst
\ebew

\begin{lemma}\label{lem_A3bis}
Let ${\cal N}$ be a von Neumann algebra with predual ${\cal N}_*$.
If a functional $\sigma$ in the unit ball of ${\cal N}_*$,
projections $r,l\in {\cal N}$
and a number $\beta\in ]0,1[$
are such that $r(\betr{\sigma})\geq 1-\beta$
and $l(\betr{\sigma^*})\geq 1-\beta$
then $\norm{\sigma -\tau}< 5\sqrt{\beta}$ where
$\tau=\frac{l\sigma r}{\norm{l\sigma r}}$.
\end{lemma}
{\em Proof}: $\norm{l\sigma r -\sigma}\leq 2\sqrt{\beta}$
by \Ref{glA3_3} and 
$\norm{\frac{l\sigma r}{\norm{l\sigma r}} -l\sigma r}
=\frac{1-\norm{l\sigma r}}{\norm{l\sigma r}}\,\norm{l\sigma r}
\leq \beta + \norm{\sigma}-\norm{l\sigma r}
\leq \beta +  \norm{\sigma-l\sigma r}
\leq \beta + 2\sqrt{\beta}$
thus $\norm{\sigma -\tau}< 5\sqrt{\beta}$.
\ebew\bigskip\\
We recall some more definitions and notation.
Let $A$ be a C$^*$-algebra. A projection $p\in A''$ is called open
if it is the limit of an increasing net of positive elements of $A$
(\cite[3.11]{Ped}, \cite{Tak}).
If $p\in A''$ is open then $B''=pA''p$ where $B=pA''p\cap A$ is
a hereditary subalgebra.
A projection $q\in A''$ is called closed if there is an open projection
$p\in A''$ such that $q=p^c$ where $p^c$ denotes the complement
$\eins-p$ of $p$. (This makes sense also if $A$ is not unital
because one always has $\eins\in A''$.)
By definition the closure $\overline{p}$ of a projection $p\in A''$
is the infimum of all closed projections majorizing $p$.
$\chi_M$ denotes the characteristic function
of a set $M$.
By functional calculus
$\chi_{]\eps,1]}(x)$ (respectively
$\chi_{[\eps,1]}(x)$) is an open (respectively closed)
projection in $A''$ if $1>\eps>0$, $x\in A$,
$0\leq x\leq\eins$, because $\chi_{]\eps,1]}$ (respectively
$\chi_{[\eps,1]}$) 
is the pointwise limit of an increasing (respectively decreasing)
sequence of continuous functions on $[0,1]$ (cf. \cite[II.3]{Ake69}).
As to $\chi_{]\eps,1]}(x)$ this is easy but as to 
$\chi_{[\eps,1]}(x)$ a bit more attention must be paid to the case
where $A$ is not unital; in this case one works with the unitisation 
$A_{\eins}$ of $A$. Therfore, in order to avoid complications in the
non-unital case \cite{AkePed}, we state Lemma \ref{lem_leins_endl}
only for unital $A$ (although the lemma holds also
in the non-unital case).

The following Lemma \ref{lem_leins_endl} is a natural generalisation
of \cite[Lem.\ 5]{Pf-V}.
\begin{lemma}\label{lem_leins_endl}
Let $A$ be a unital \Cstalg.
For each $\eps>0$ and each $n\in\N$ there is $\delta=\delta(n,\eps)>0$
with the following property.

If there are functionals
$\phi_1,\ldots,\phi_n$ in the unit ball of $A'$ and open projections
$s,t\in A''$ such that
\bgl
(1-\delta)\sum_1^n\betr{\alpha_k} \leq
\norm{\sum_1^n \alpha_k\, t\phi_k s} \leq \sum_1^n \betr{\alpha_k}
\falle (\alpha_{k})\aus {\C}
\label{gl_leins_endl2}
\egl
then there are open projections $p_1,\ldots,p_n\in sA''s$ with pairwise
orthogonal closures in $sA''s$ 
and open projections $q_1,\ldots,q_n\in tA''t$
with pairwise orthogonal closures in $tA''t$ such that
\bgl
p_k(\betr{\phi_k})      &>&      (1-\eps)\norm{\phi_k}
\label{gl_leins_endl5}\\
q_k(\betr{\phi_k^*})      &>&      (1-\eps)\norm{\phi_k}
\label{gl_leins_endl6}
\egl
for $k=1,\ldots,n$.\\
In particular the $\phi_k$ are close to normalized orthogonal elements
$\psi_k$ on $sAt$ in the sense that $\norm{\phi_k -\psi_k}<5\sqrt{\eps}$
where $\psi_k={q_k \phi_k p_k}/{\norm{q_k \phi_k p_k}}$ are
normalized and pairwise orthogonal with left (right) supports
majorized by $t$ (by $s$).
\end{lemma}
{\em Proof}:
(a)
First we suppose $s=t=\eins$ and deal only with the special case\\
(a1) of positive functionals $\phi_k$.\\
Let $\eps>0$.
For $n=1$ choose an $x\geq 0$ in the unit ball of $A$ such that
$\phi_1(x) > 1-\eps$ and set $p_1=q_1=\chi_{]0,1]}(x)$,
$\delta(1,\eps)=\eps$.

Suppose now that the assertion holds true (for positive functionals, for
$s=t=\eins$ and)
for $n\in\N$. By hypothesis on $n$ we choose $\delta_{n}=\delta(n,\eps)$
We define $\delta_{n+1}>0$ such that
$\delta_{n+1} +(32n\delta_{n+1})^{1/2} <\delta_{n}$.
Consider positive functionals $\phi_k$, $k=1,\ldots,n+1$,
in the unit ball of $A$ such that
$\norm{\sum_{1}^{n+1} \alpha_k \phi_k} \geq (1-\delta_{n+1}) 
\sum_{1}^{n+1} \betr{\alpha_k}$.
Set 
$\sigma=\frac{1}{n}\sum_1^n\phi_k$ and $\tau=\phi_{n+1}$.
Then $(1-\delta_{n+1})(\betr{\alpha}+\betr{\beta})
\leq
\norm{\alpha\sigma+\beta\tau}
\leq \betr{\alpha}+\betr{\beta}$
for all scalars $\alpha,\beta$.
In particular $\norm{\sigma-\tau}\geq 2(1-\delta_{n+1})$.
There is a selfadjoint normalized element $x\in A$ such that 
$(\sigma-\tau)(x)>2(1-2\delta_{n+1})$.
Decompose $x=x^+ +x^-$ in its negative and positive parts.
Then
$(\sigma-\tau)(x)=(\sigma(x^+)+\tau(x^-))-(\sigma(x^-)+\tau(x^+))
>2(1-2\delta_{n+1})$ whence
$\sigma(x^+)>2(1-2\delta_{n+1})-\tau(x^-)
>1-4\delta_{n+1}$ and similarly
$\phi_{n+1}(x^-)>1-4\delta_{n+1}$.
Together with $\phi_k(x^+)\leq 1$ this gives $\phi_k(x^+)>1-4n\delta_{n+1}$
for all $k=1,\ldots,n$ because otherwise one would have 
$n\sigma(x^+)\leq 1-4n\delta_{n+1} +(n-1)= n(1-4\delta_{n+1})$.
Define $p=\chi_{]\eta,1]}(x^+)$, $p_{n+1}=\chi_{]\eta,1]}(x^-)$
where $\eta>0$ is such that
\bgl
p(\phi_k)&>&1-4n\delta_{n+1}
\,\,\,\,\,\,\,\,\,\,\mbox{ for } k=1,\ldots, n
\label{gl_leins_endl4}\\
\mbox{and }\;\;\;\;\;\;\;
p_{n+1}(\phi_{n+1})&>&1-4\delta_{n+1}>1-\eps \nonumber
\egl
By functional calculus the projections $p$ and $p_{n+1}$
are open and have orthogonal closures.
$B=pA''p\cap A \aus A$ is a hereditary subalgebra of $A$.
This explains the equality sign in the following formula:
\bglst
\norm{\sum_{1}^{n}\alpha_k \phi_k\eing{B}}_B
&=& 
\norm{\sum_{1}^{n}\alpha_k(p \phi_k p)}
\geq
\norm{\sum_{1}^{n}\alpha_k \phi_k}
-\norm{\sum_{1}^{n}\alpha_k (\phi_k - p \phi_k p)}\\
&\stackrel{\Ref{glA4_3}\Ref{gl_leins_endl4}}{>}&
(1-\delta_{n+1})\sum_{1}^{n}\betr{\alpha_k}
      -\sqrt{32n\delta_{n+1}}\sum_{1}^{n}\betr{\alpha_k}\\
&>&(1-\delta_{n})\sum_{1}^{n}\betr{\alpha_k}.
\eglst
By induction hypothesis applied to $B$ and to $\phi_{k}\eing{B}$
one gets $n$ open
projections $p_1,\ldots,p_n\in B''$ with pairwise orthogonal closures
in  $B''$ - whence in $A''$ -
such that \Ref{gl_leins_endl5} holds for $k=1,\ldots,n+1$.
Furthermore, \Ref{gl_leins_endl6} holds with
$q_k=p_k$ because we have supposed $\phi_k\geq 0$.
This proves the case where $s=t=\eins$ for positive functionals $\phi_k$.\\
(a2) For the case of arbitrary functionals (but still with $s=t=\eins$)
suppose that the lemma is false. Then there are $\eps>0$,
a sequence $(A_i)$ of \Cstalg s, and $\phi_{k,i}\in (A_i)_{\eins}$
such that for each $i\in\N$, 
\bgl
\Bigl(1-\frac{1}{i}\Bigr)\sum_{k=1}^n \betr{\alpha_k}
<
\norm{\sum_{k=1}^n \alpha_k \phi_{k,i}}
 \leq \sum_{k=1}^n \betr{\alpha_k}
\,\,\,\,\falle (\alpha_k)\subset\in\C
                            \label{gl_leins_endl7}
\egl
but for each $i\in\N$ the $\phi_{k,i}$ are far from orthogonal
functionals, more precisely
\bgl
\min_{k\leq n}p_{k,i}(\betr{\phi_{k,i}}) \leq 1-\eps
\,\,\,\,\,\mbox{or}\,\,\,\,\,
\min_{k\leq n}q_{k,i}(\betr{\phi_{k,i}^*}) \leq 1-\eps
\label{gl_leins_endl1}
\egl
for all sequences $(p_{k,i})_{k=1}^n$ and $(q_{k,i})_{k=1}^n$ of
open projections with orthogonal closures in $A_i''$.

We recall some basic facts on ultraproducts (see e.g. \cite{Hei1}).
If ${\cal U}$ is an
ultrafilter on an index set $I$ the ultraproduct $X=(X_i)/{{\cal U}}$
of a family $(X_i)_{i\in I}$ of Banach spaces is defined as the
quotient $\lunendl(X_i)/c_0(X_i)$ where 
$\lunendl(X_i)=\{(x_i)_{i\in I}\mdE
\norm{(x_i)};=\sup_{\cal U}\norm{x_i}<\infty\}$
and
$c_0(X_i)=\{(x_i)_{i\in I}\in \lunendl(X_i)\mdE \lim_{{\cal U}}\norm{x_i}=0\}$.
With the quotient norm $X$ becomes a Banach space. By $[x_i]_{{\cal U}}$ we denote the equivalence
class represented by $(x_i)_{i\in I}\in \lunendl(X_i)$.
One has $\norm{[x_i]_{{\cal U}}}=\lim_{{\cal U}}\norm{x_i}$ independently of the
representative of $[x_i]_{{\cal U}}$.
The ultraproduct $(X_i')/{{\cal U}}$ of the duals can be identified isometrically
with a closed subspace of the dual $X'$ via
$[x_i']_{{\cal U}}([x_i]_{{\cal U}})=\lim_{{\cal U}}x_i'(x_i)$.
An ultraproduct $A=(A_i)/{{\cal U}}$ of a family of \Cstalg s $A_i$ is canonically
a \Cstalg\ with pointwise multiplication and involution because
in this case the null space $c_0(X_i)$ is an ideal in $A$.

Let now $I=\N$ and ${\cal U}$ be a free ultrafilter on $\N$
and set $A=(A_i)/{{\cal U}}$. For each element $\psi\in A'$ of the form
$\psi=[\psi_{i}]_{{\cal U}}$ we have $\betr{\psi}=[\betr{\psi_{i}}]_{{\cal U}}$
and $\betr{\psi^*}=[\betr{\psi_{i}^{*}}]_{{\cal U}}$.
[To see this choose $a=[a_{i}]_{{\cal U}}$ in the unit ball
of $A$ such that $\norm{a_{i}}=1$ and
$\psi(a)=\lim_{{\cal U}} \psi_{i}(a_{i})=\lim_{{\cal U}} \norm{\psi_{i}}=1$.
Then $\betr{\psi}=a\psi=[a_{i}\psi]_{\cal U}$ and
$\norm{\,\betr{\psi_{i}}-a_{i} \psi_{i}}
     \leq (2\betr{\,\norm{\psi_{i}}-\psi_{i}(a_{i})})^{1/2} \gen 0$
by \Ref{glA3_2} of Lemma \ref{lem_A3}
hence $\betr{\psi}=[\betr{\psi_{i}}]_{{\cal U}}$.
For $\betr{\psi^*}=[\betr{\psi_{i}^{*}}]_{{\cal U}}$ the proof is
analogous.]\\
The $n$ functionals $\betr{\phi_{k}}=[\betr{\phi_{k,i}}]_{{\cal U}}$
and the $n$ functionals $\betr{\phi_{k}^*}=[\betr{\phi_{k,i}^{*}}]_{{\cal U}}$
are pairwise orthogonal because the $\phi_{k}=[\phi_{k,i}]_{{\cal U}}$ are so
by \Ref{gl_leins_endl7}.
Since the statement of the Lemma has been proved for
positive functionals (and for $s=t=\eins$)
there are, for each $i\in I$, two finite sequences
$(p_{k,i})_{k=1}^{n}$, $(q_{k,i})_{k=1}^{n}$ of open projections with pairwise
orthogonal closures in $A_i''$ such that 
$\lim_{{\cal U}}p_{k,i}(\phi_{k,i}) = [p_{k,i}]_{{\cal U}}(\phi_{k})=1$ and
$\lim_{{\cal U}}q_{k,i}(\phi_{k,i}^{*}) =
[q_{k,i}]_{{\cal U}}(\phi_{k}^{*})=1$ for $k=1,\ldots,n$
which contradicts \Ref{gl_leins_endl1} and thus proves the lemma 
for the case where $s=t=\eins$.\\
(b) Now we turn to the general case of arbitrary open  projections
$s,t\in A''$.
We assume without loss of generality that $\norm{\phi_k}=1$.
After what has been proved in (a), we further assume without
loss of generality that the $\phi_k$ are pairwise orthogonal
because \Ref{gl_leins_endl2} remains valid if $t\phi_k s$ is replaced
by $\phi_k$.
Then the $\betr{\phi_k}$ are orthogonal
and so are the $\betr{\phi_k^*}$.
Thus we may further assume that
\bgl
\norm{\sum_1^n \alpha_k\betr{\phi_k}\,}
=\sum_1^n \betr{\alpha_k}
\,\,\,\,\,\,\mbox{and}\,\,\,\,\,\,
\norm{\sum_1^n \alpha_k\betr{\phi_k^*}\,}
=\sum_1^n \betr{\alpha_k}
\label{gl_leins_endl2a}
\egl
for all scalars $\alpha_k\in\C$.
Let $\phi_k=u_k\betr{\phi_k}$ be the polar decomposition then
\bglst
1-\delta
\stackrel{\Ref{gl_leins_endl2}}{\leq}
 \norm{t\phi_k s}
=\norm{t u_k \betr{\phi_k} s}
=\norm{\,\betr{tu_k\betr{\phi_k} s}\,}
=s(\betr{(tu_k\betr{\phi_k} s)})
\leq
s(\betr{\phi_k})
\eglst
where the last inequality follows from \cite[III.4.9]{Tak}.
Analogously $t(\betr{\phi_k^*})\geq 1-\delta$. Hence by \Ref{glA3_3}
of Lemma \ref{lem_A3}, if $\delta$ is small enough,
$t\phi_k s$ is a small perturbation of $\phi_k$ and since the absolut
value is norm continuous  on von Neumann preduals \cite[III.4.10]{Tak},
$\betr{t\phi_k s}$ is a small perturbation of $\betr{\phi_k}$, and
$\betr{s\phi_k^* t}$ is a small perturbation of $\betr{\phi_k^*}$.
Finally, in view of \Ref{gl_leins_endl2a}  we may
without loss of generality replace \Ref{gl_leins_endl2} by
\bgl
(1-\delta)\sum_1^n\betr{\alpha_k} \leq
\norm{\sum_1^n \alpha_k\, \betr{t\phi_k s}\,}
\leq \sum_1^n \betr{\alpha_k}
\falle (\alpha_{k})\aus {\C}
\label{gl_leins_endl2b}\\
(1-\delta)\sum_1^n\betr{\alpha_k} \leq
\norm{\sum_1^n \alpha_k\, \betr{s\phi_k^* t}\,}
\leq \sum_1^n \betr{\alpha_k}
\falle (\alpha_{k})\aus {\C}
\label{gl_leins_endl2c}
\egl
in the statement of the lemma. It remains to apply part (a) to the
hereditary subalgebra $sA''s\cap A$ and to $tA''t\cap A$
(because the support projection of $\betr{t\phi_k s}$
(of $\betr{s\phi_k^* t}$) is majorized by $s$ (by $t$)). This yields
the desired open projections $p_k\in (sA''s\cap A)''=sA''s$
and $q_k\in tA''t$ satisfying
\Ref{gl_leins_endl5} and \Ref{gl_leins_endl6}
if $\delta$ is small enough.\\
The last assertion of the lemma is immediate from Lemma \ref{lem_A3bis}.
\ebew\bigskip\\
\begin{coro}\label{coro_Wst}
Let ${\cal N}$ be a von Neumann algebra.
For each $\eps>0$ and each $n\in\N$ there is $\delta=\delta(n,\eps)>0$
with the following property.

If there are functionals
$\phi_1,\ldots,\phi_n$ in the unit ball of ${\cal N}_*$
and (arbitrary) projections $s,t\in {\cal N}$ such that
\bgl
(1-\delta)\sum_1^n\betr{\alpha_k} \leq
\norm{\sum_1^n \alpha_k\, t\phi_k s} \leq \sum_1^n \betr{\alpha_k}
\falle (\alpha_{k})\aus {\C}
\label{gl_leins_endl2bis}
\egl
then there are pairwise orthogonal projections
$p_1,\ldots,p_n\in s{\cal N}s$ 
and pairwise orthogonal projections $q_1,\ldots,q_n\in t{\cal N}t$ such that
\bgl
\norm{\phi_k -\psi_k}<\eps&&\mbox{for }k=1,\ldots,n
\egl
where $\psi_k={q_k \phi_k p_k}/{\norm{q_k \phi_k p_k}}$.
\end{coro}
{\em Proof}:
For $s=t=\eins$ the assertion is immediate from Lemma \ref{lem_leins_endl}
and Lemma \ref{lem_A3bis}.
For arbitrary projections $s,t\in{\cal N}$ we proceede  as in part (b)
of the proof of Lemma \ref{lem_leins_endl} in order to show that
\Ref{gl_leins_endl2bis} can be replaced by \Ref{gl_leins_endl2b}
and \Ref{gl_leins_endl2c} and
to apply this to the subalgebras $s{\cal N}s$ and $t{\cal N}t$.
\ebew\bigskip\bigskip\bigskip\\
{\bf\S 4 Proof of Theorem \ref{prop_leins_asy}}
\bigskip\\
Without loss of generality we assume that $\norm{\phi_m}=1$ for all
$m\in\N$.
Let $(\eta_n)$ be a sequence of positive numbers such that $\sum\eta_n$
converges.

By induction on $n=1,2,\ldots$ we construct
an increasing sequence $(m_n)$ in $\N$,
functionals $\psi_{m_{k}}^{(n)}\in{\cal N}_*$ for $k=1,\ldots,n$,
such that for all $n\in\N$:
\bgl
\betr{\psi_{m_{k}}^{(n)}} &\bot& \betr{\psi_{m_{l}}^{(n)}}
           \,\,\,\,\,\,\,   k,l=1,\ldots,n,\,\,\,k\neq l,   
                         \label{glprop_leins_asy1}\\
\norm{\psi_{m_{k}}^{(n)}}  &=& 1  
  \,\,\,\,\,\,\,\,\,\,\,\,\,\,\,\,\,   k=1,\ldots,n
                      \label{glprop_leins_asy3}\\
\norm{\psi_{m_{k}}^{(n)} - \psi_{m_{k}}^{(n-1)}} &<&
                     \eta_n 
\,\,\,\,\,\,\,\,\,\,\,\,\,\, k=1,\ldots,n-1,  
                         \label{glprop_leins_asy2}\\
\norm{\psi_{m_{n}}^{(n)} - \phi_{m_{n}}} &<&     \eta_n
                     \label{glprop_leins_asy4}.
\egl
For $n=1$ one may simply set $\psi_{m_{1}}^{(1)}=\phi_1$;
(\ref{glprop_leins_asy1}, $n=1$) and
(\ref{glprop_leins_asy2}, $n=1$) are void,
(\ref{glprop_leins_asy3}, $n=1$) and (\ref{glprop_leins_asy4}, $n=1$)
are trivial.\\
Induction step $n\mapsto n+1$.\\
Suppose $m_k$ and $\psi_{m_{k}}^{(n)}$ to be constructed
for $k=1,\ldots,n$  according to
\Ref{glprop_leins_asy1} - \Ref{glprop_leins_asy4}.\\
Choose $\delta_1=\delta(n,\eta_{n+1}/2)>0$ according to
Corollary \ref{coro_Wst}
such that furthermore $\delta_1<\eta_{n+1}/2$.
Let $j\in\N$ be such that $(2/j)^{1/2}<\delta_1$.
Now, again according to  Corollary \ref{coro_Wst}, choose
$\delta_0=\delta(nj,\eta_{n+1})$.\\
Since $(\phi_m)$ spans $\leins$ almost isometrically
there is an index $m_0$ such that $(\phi_m)_{m\geq m_0}$ spans $\leins$
$(1-\delta_0)$-isomorphically.
By Corollary \ref{coro_Wst} (with $s=t=\eins$, $\delta=\delta_0$)
we find a finite set $N\subset\N$ of
cardinality $nj$ (for example $N=\{m_0 +1,\ldots,m_0 +nj\}$),
a finite sequence of orthogonal projections
$(p_m)_{m\in N}$ in ${\cal N}$ such that
\bgl
\norm{\phi_m - \frac{\phi_m p_m}{\norm{\phi_m p_m}}} 
  &<&  \eta_{n+1}
             \,\,\,\,\falle m\in N_{n+1}.   \label{glprop_leins_asy5}
\egl
Set $\phi=\sum_{k=1}^{n} \betr{\psi_{m_{k}}^{(n)}}$.
$\phi$ is positive. We have
$(\sum_{m\in N}{p}_m)(\phi)\leq\norm{\phi}\leq n$.
Thus there is an index $m_{n+1}\in N$ such that 
$0\leq {p}_{m_{n+1}}(\phi) \leq 1/j$ and
$0\leq {p}_{m_{n+1}}(\betr{\psi_{m_{k}}^{(n)}}) \leq 1/j$ for $k=1,\ldots,n$.
We set $s=\eins - \overline{p}_{m_{n+1}}$ and define
$\tilde{\psi}_{m_{k}}^{(n+1)}=\psi_{m_{k}}^{(n)} s$
for $k=1,\ldots,n$ and
$$\psi_{m_{n+1}}^{(n+1)}
=\frac{\phi_{m_{n+1}}p_{m_{n+1}}}{\norm{\phi_{m_{n+1}}p_{m_{n+1}}}}.$$
Then (\ref{glprop_leins_asy3}, $n+1$) holds for $k=n+1$ and
(\ref{glprop_leins_asy4}, $n+1$) holds by \Ref{glprop_leins_asy5}.
We have
$s(\betr{\psi_{m_k}^{(n)}})=\norm{\psi_{m_k}^{(n)}}
-{p}_{m_{n+1}}(\betr{\psi_{m_k}^{(n)}})
\geq 1-1/j$ by \Ref{glprop_leins_asy3}.
From this and \Ref{glA3_3} one gets that
\bgl
\norm{\tilde{\psi}_{m_{k}}^{(n+1)} -\psi_{m_{k}}^{(n)}}
\leq (2/j)^{1/2} <\delta_1
   <\frac{\eta_{n+1}}{2},\,\, \,\,\,k=1,\ldots,n.
\label{glprop_leins_asy6}
\egl
Thus, up to $\delta_1$  the $\tilde{\psi}_{m_{k}}^{(n+1)}$ are near to an
isometric copy of $\leins_n$ because
\bglst
\sum_{k=1}^n \betr{\alpha_k}
&\geq&\norm{\sum_{k=1}^n \alpha_k \tilde{\psi}_{m_{k}}^{(n+1)}}
=\norm{\sum_{k=1}^n \alpha_k \tilde{\psi}_{m_{k}}^{(n+1)}s}\\
&\geq& \norm{\sum_{k=1}^n \alpha_k \psi_{m_{k}}^{(n)}} - 
                \norm{\sum_{k=1}^n \alpha_k(\tilde{\psi}_{m_{k}}^{(n+1)} 
                                         -\psi_{m_{k}}^{(n)})}\\
&\stackrel{\Ref{glprop_leins_asy6}}{\geq}&
\norm{\sum_{k=1}^n \alpha_k \psi_{m_{k}}^{(n)}} 
          -(2/j)^{1/2}\sum_{k=1}^n \betr{\alpha_k}\\
&\stackrel{(\ref{glprop_leins_asy1}, \ref{glprop_leins_asy3})}{=}&
\bigl( 1-(2/j)^{1/2}\bigr)\sum_{k=1}^n \betr{\alpha_k}\\
&>& (1-\delta_1)\sum_{k=1}^n \betr{\alpha_k}.
\eglst
It remains to apply Corollary \ref{coro_Wst} another time
(with $t=\eins$, $\delta=\delta_1$) in order
to get small normalized orthogonal perturbations
$\psi_{m_{k}}^{(n+1)}$ - whence (\ref{glprop_leins_asy3}, $n+1$)
for $k\leq n$ -
of the $\tilde{\psi}_{m_{k}}^{(n+1)}$
whose right supports are majorized by $s$ and thus
orthogonal to the right support of $\psi_{m_{k}}^{(n+1)}$ such that
$\norm{\psi_{m_{k}}^{(n+1)}-\tilde{\psi}_{m_{k}}^{(n+1)}}<\eta_{n+1}/2$ for $k=1,\ldots,n$.
Together with \Ref{glprop_leins_asy6} this gives (\ref{glprop_leins_asy2}, $n+1$). Finally
one verifies (\ref{glprop_leins_asy1}, $n+1$) by observing that the 
support projections of the $\betr{\psi_{m_{k}}^{(n+1)}}$ are the right supports of
the $\psi_{m_{k}}^{(n+1)}$. This ends the induction.

By construction, $(\psi_{m_{k}}^{(n)})_{n\in\N}$ is a Cauchy sequence for each
$k$ because $\norm{\psi_{m_{k}}^{(n)}-\psi_{m_{k}}^{(i)}}\leq \sum_{l=i+1}^{n}\eta_l\gen 0$
as $n>i\gen\infty$. Let $\psi_k=\lim_{n}\psi_{m_{k}}^{(n)}$
be its limit. Then
$\norm{\psi_k-\phi_{m_k}}\leq\norm{\phi_{m_k}-\psi_{m_{k}}^{(k)}}+
                  \norm{\psi_{m_{k}}^{(k)}-\lim_{n}\psi_{m_{k}}^{(n)}}
         \leq \eta_k +\sum_{l=k+1}^{\infty}\eta_l\gen 0$
as $k\gen\infty$.
The $\psi_k$ have pairwise orthogonal right supports because
by continuity of the absolute value (\cite[III.4.10]{Tak}),
if $k\neq l$ one has 
\bglst
\norm{\,\betr{\psi_k} - \betr{\psi_l}\,} 
&=&
 \lim_{n\gen\infty}\norm{\,\betr{\psi_{m_{k}}^{(n)}}-\betr{\psi_{m_{l}}^{(n)}}\,}\\
&\stackrel{\Ref{glprop_leins_asy1}}{=}&
\lim_{n\gen\infty}\norm{\,\betr{\psi_{m_{k}}^{(n)}}\,}
+\norm{\,\betr{\psi_{m_{l}}^{(n)}}\,}
=\norm{\psi_k}+\norm{\psi_l}.
\eglst
So far we have proved that if $(\phi_m)$ spans $\leins$
almost isometrically then there is a subsequence $(\phi_{m_k})$ and
there are pairwise orthogonal projections $s_k\in {\cal N}$
(namely the right support projections of the $\psi_k$) such that
$\norm{\phi_{m_k}-\phi_{m_k}s_k}
\leq \norm{\phi_{m_k}-\psi_k}+\norm{\psi_k s_k - \phi_{m_k}s_k}
\leq 2\norm{\phi_{m_k}-\psi_k}\gen 0$.
Since $(\phi_{m_k}^*)$ spans $\leins$ almost isometrically, too, there
are pairwise orthogonal projections $t_l\in {\cal N}$ such that
$\norm{\phi_{m_{k_l}}^* -\phi_{m_{k_l}}^* t_l}\gen 0$
for an appropriate sequence $(m_{k_l})$ in $\N$. Set
$\tilde{\phi}_l = t_l \phi_{m_{k_l}} s_{k_l}$.
Then $\norm{\phi_{m_{k_l}}-\tilde{\phi}_l}
\leq \norm{\phi_{m_{k_l}}-\phi_{m_{k_l}} s_{k_l}}
   + \norm{(\phi_{m_{k_l}}-t_l \phi_{m_{k_l}}) s_{k_l}}
\leq \norm{\phi_{m_{k_l}}-\phi_{m_{k_l}} s_{k_l}}
    +\norm{\phi_{m_{k_l}}^* -\phi_{m_{k_l}}^* t_l}
\gen0$.\\
The second statement of the theorem is trivial by the definiton
of the $\tilde{\phi}_l$.
This ends the proof.
\ebew
\bigskip\\
From Remark 2 after the proof of Theorem \ref{theo_L1} at the end
of the next section it follows that
Theorem \ref{prop_leins_asy} does not hold for unbounded sequences
$(\phi_m)$.
\bigskip\bigskip\bigskip\\
{\bf\S 5 Proof of Theorem \ref{theo_L1}}
\bigskip\\
(i) $\folgt$ (ii). Let $(x_{n_{k}})$ be a subsequence of $(x_n)$.
If $(x_{n_{k}})$ contains
a sequence $(x_{n_{k_l}})$ such that $x_{n_{k_l}}=0$ for all $l\in\N$
then we simply choose $y_l=0$ for $l\in\N$.
Otherwise we may (pass to another subsequence and) suppose that
$\norm{x_{n_{k}}}_1\neq0$ for all $k\in\N$.
By norm density of $\Leins \cap \Lunendl$ in $\Leins$ and the fact
that the norm topology is finer than the measure topology we may suppose
without loss of generality that
$0\neq\norm{x_{n_{k}}}_{\infty}<\infty$ for all $k\in\N$.
We set $\eps_l=2^{-l}/\tau(\eins)$ for $l\in\N$.

By induction over $l\in\N$ we construct
a strictly increasing subsequence $(n_{k_l})$ of $(n_k)$,
projections $p_l\in {\cal N}$
and positive numbers $\delta_l$
such that for all $l\in\N$
\bgl
\tau(p_l)<\delta_l 
\,\,\,\,\,\mbox{where} \,\,\,\,
p_l=\chi_{]\eps_l,\infty[}(\betr{x_{n_{k_l}}})
                          \label{glvN8}
\egl
and where
\bgl
\delta_l &=& \frac{2^{-l}}{\max_{\,1\leq m\leq l-1}\norm{x_{n_{k_m}}}_{\infty}},
\,\,\,\,\mbox{if}\,\,\, l\geq 2.
                                              \label{glvN9}
\egl
For $l=1$ we choose $n_{k_1}=n_1$ and any $\delta_1>\tau(p_1)$.
For the induction step $l\mapsto l+1$
we suppose $n_{k_m}$, $p_{m}$, and $\delta_m$
to be constructed for $m=1,\ldots,l$,
we define $\delta_{l+1}$ by \Ref{glvN9}
and choose $n_{k_{l+1}}$ such that
$$\tau(\chi_{]\eps_{l+1},\infty[}(\betr{x_{n_{k_{l+1}}}}) <\delta_{l+1}$$
which is possible because $x_n\gentau0$.
We define $p_{l+1}$ by \Ref{glvN8}.
This settles (\ref{glvN8}, $l+1$) and ends the induction.\\
By \Ref{glvN9} we have
\bglst
\delta_{l+1+r}= 
\frac{2^{-(l+1+r)}}{\max_{\,m\leq l+r}\norm{x_{n_{k_m}}}_{\infty}}
\leq 2^{-(r+1)}\frac{2^{-l}}{\norm{x_{n_{k_l}}}_{\infty}}
\eglst
for $r\in\N\cup\{0\}$
which gives
\bgl
\sum_{m\geq l+1} \delta_m 
=
\sum_{r\geq 0}\delta_{l+1+r}
\leq\frac{2^{-l}}{\norm{x_{n_{k_l}}}_{\infty}}.
\label{glvN1}
\egl                                                              
Put
$q_l=1-\bigvee_{m\geq l+1}p_m$ and
$\tilde{y}_l=x_{n_{k_l}}\bigl(p_l\wedge q_l)$.
By construction the $\tilde{y}_l$ have pairwise orthogonal right support projections
and their left support projections are majorized by the ones of the $x_{n_{k_l}}$.
We show that $\norm{x_{n_{k_l}}-\tilde{y}_l}_1\gen 0$.
In order to save indices we use the abbreviations
$x=x_{n_{k_l}}$, $p=p_l$, $q=q_l$, $\tilde{y}=\tilde{y}_l$
until the end of formula \Ref{glvN15}:
\bgl
\norm{x-\tilde{y}}_1
&\leq& \norm{x-x p}_1 +\norm{x p -\tilde{y}}_1
                \nonumber\\
&=&\norm{x(\eins-p)}_1 +\norm{x(p-(p\wedge q))}_1
                \nonumber\\
&\leq& \norm{x(\eins-p)}_{\infty}\,\tau(\eins)
              +\norm{x}_{\infty}\,\tau(p-(p\wedge q))
                \nonumber\\
&\stackrel{(*)}{=}&
   \norm{\,\betr{x}\,\chi_{[0,\eps_l]}(\betr{x})}_{\infty}\,\tau(\eins)
             +\norm{x}_{\infty}\,\tau\Bigl((p\vee q)-q\Bigr)
                \nonumber\\
&\leq& \eps_l \,\tau(\eins) 
           + \norm{x}_{\infty}\,\tau(\eins-q)
                \nonumber\\
&\leq& \eps_l \,\tau(\eins) 
           + \norm{x}_{\infty}\,\Bigl(\sum_{m\geq l+1}\tau(p_m)\Bigr)
                \nonumber\\
&\stackrel{\Ref{glvN8}, \Ref{glvN1}}{\leq}&
2^{-(l-1)}.
\label{glvN15}
\egl
For $(*)$ we used that $p-(p\wedge q)$ and $(p\vee q)-q$ are equivalent
projections for any two projections $p,q$ (\cite[V.1.6]{Tak}) hence
$\tau(p-(p\wedge q))=\tau((p\vee q)-q)$.
\smallskip

So far we have proved that given a $\tau$-null subsequence
$(x_{n_k})$ there are $x_{n_{k_l}}$ and there are $\tilde{y}_l$
whose right supports are orthogonal
and whose left supports are majorized by the left supports of the
$x_{n_{k_l}}$ such that $\norm{x_{n_{k_l}}-\tilde{y}_l}_1\gen0$.
In particular, $\tilde{y}_l\gentau0$ whence
$\tilde{y}_l^*\gentau0$. Thus we can apply the same reasoning
(up to passing to appropriate subsequences) in order to find 
perturbations $y_l^*$ of the $\tilde{y}_l^*$ which have both orthogonal right
and orthogonal left supports such that
$\norm{\tilde{y}_l-y_l}_1=\norm{\tilde{y}_l^*-y_l^*}_1\gen0$
hence $\norm{x_{n_{k_l}}-y}_1\gen0$.
This ends the proof of (i) $\folgt$ (ii).\\
(ii) $\folgt$ (i) Since $\tau$ is finite
and the $y_l$ are pairwise orthogonal we have that $y_l\gentau0$.
And $\norm{x_{n_{k_l}}-y_l}\gen0$ entails $x_{n_{k_l}}-y_l\gentau0$
hence $x_{n_{k_l}}\gentau0$.
Thus each subsequence of $(x_n)$ contains a subsequence which converges
to $0$ in measure whence $x_n\gentau0$.\\
(ii) $\folgt$ (iii) follows from Lemma \ref{lem_asy_stabil}:
Suppose (ii) holds and $\inf\norm{x_{n_k}}_1 >0$ for a subsequence
$(x_{n_k})$ of $(x_n)$. Then by (ii), there are orthogonal $y_l$
and there is $(x_{n_{k_l}})$ such
that $\norm{x_{n_{k_l}}-y_l}_1\gen 0$. One may suppose that
$\inf\norm{y_l}_1 >0$ hence $(y_l)$ spans $\leins$
isometrically. Thus by Lemma \ref{lem_asy_stabil}, the sequence
$(x_{n_{k_l}})=(y_l+(x_{n_{k_l}}-y_l))$ spans
$\leins$ almost isometrically.\\
(iii) $\folgt$ (iv): Von Neumann preduals are L-embedded spaces
\cite[IV.1.1]{HWW},
thus by \cite{Pfi-Fixp} each sequence spanning $\leins$ almost
isometrically admits a subsequence spanning $\leins$ asymptotically.\\
(iv) $\folgt$ (iii) is trivial.\\
(iii) $\folgt$ (ii) follows from Theorem \ref{prop_leins_asy}.\smallskip\\
See the following Remark 2 for an example which shows
that in general (iii) does not imply (i), (ii)
for unbounded sequences $(x_n)$.
\ebew\bigskip\medskip\\
Remarks:\\
1. As an illustration of how to get an orthogonal subsequence
consider the sequence $x_n=n\,\charF{[0,\,1/n]}$ in $\Leins([0,1])$.
One may take, for example,
$y_l=x_{n_l} \eins_{]1/n_{l+1},1]}=n_l  \eins_{]1/n_{l+1},1/n_l]}$
where $n_l=2^{(2^{l})}$.\medskip\\
2. In general (iii) does not imply (i), neither (ii),
if the sequence $(x_n)$ is unbounded.
Take the bounded sequence $x_n=n^2 \charF{[1/n+1,\;1/n[}+ \frac{1}{n}$ in $\Leins([0,1])$.
It converges to zero in measure and does not contain a norm null sequence.
Hence by (i)$\folgt$(iii) an appropriate subsequence $({x_{n_{k}}})$
spans $\leins$ almost isometrically.
Thus the unbounded sequence $(n_k^2 x_{n_{k}})$ satisfies (iii)
but not (i). It cannot satisfy (ii) either because
(ii) $\gdw$ (i) holds also for unbounded sequences.
This means in particular that Theorem \ref{prop_leins_asy}
does not hold for unbounded sequences $(\phi_m)$.\\
3. A few straightforward modifications show that (i)$\gdw$(ii)
holds accordingly also for $\LL^p({\cal N},\tau)$, $1\leq p<\infty$.
(Cf.\ \cite{Suk}.)
\newpage\noindent
{\bf\S 6 $\leins$-copies in the dual of C$^*$-algebras,
proof of Proposition \ref{theo_A1}}
\bigskip\bigskip\\
The proof of the main result of \cite{Pf-V} gives the following:
Let $(\phi_m)\subset A'$ be a bounded sequence of selfadjoint
functionals on a \Cstalg\ $A$, let $\eps>0$. If  $(\phi_m)$ spans
$\leins$ $r$-isomorphically ($0<r<1$) then there is a subsequence 
$(\phi_{m_n})$ and there is a sequence $(x_n)$ of pairwise orthogonal
normalized selfadjoint elements of $A$ such that
$\phi_{m_n}(x_n) > (1-\eps)r \norm{\phi_{m_n}}$.
This amounts to saying that
$\betr{\phi_{m_n}}(\betr{x_n})> (1-\eps)r \norm{\phi_{m_n}}$
(to see this it is enough to decompose both $\phi_{m_n}$ and $x_n$
in their positive and negative parts)
or that, via Lemma \ref{lem_A3},
$\norm{\phi_{m_n}-a_n \phi_{m_n} a_n}\gen0$ where $a_n=\betr{x_n}$.
This is Lemma \ref{lemma_A1} for selfadjoint $\phi_m$ with the better
factor $r$ instead of $r^2$ in $\Ref{gl_leins_r2}$ and
$\Ref{gl_leins_r3}$.

With Lemma \ref{lem_leins_endl} at one's disposal,
the proof of Lemma \ref{lemma_A1}
- and thus of Proposition \ref{theo_A1} -
is a straightforward
modification of \cite{Pf-V} and
gives a kind of quantitative version of \cite{Pf-V}
which holds for arbitrary functionals, not only selfadjoint ones.
(We give the entire proof of Proposition \ref{theo_A1}
not only for the sake of completeness but
also because it is quite lengthy wherefore the usual argument
"The details are left to the reader" would be exaggerated.)
Yet, it does not complete the subject "perturbations of
$\leins$-copies in C$^*$-algebras"
as at least two questions remain open.

Firstly, is it necessary in Theorem \ref{prop_leins_asy} or in Proposition
\ref{theo_A1} to pass to subsequences?
In the commutative case it is not, as a
result of Dor \cite{Dor} shows that, if 
${\cal N}_*=\Leins([0,1])$ contains a $(1-\delta)$-isomorphic copy
of $\leins$ then the whole canonical basis of this copy
can be perturbed in norm so to span $\leins$ isometrically
with the perturbation smaller than
$\delta'$ and $\delta'\gen0$ as $\delta\gen0$.
Furthermore Arazy \cite{Araz-embed}
proved that if the predual of an arbitrary
von Neumann algebra ${\cal N}$ contains a $(1-\delta)$-copy
of $\leins$ then the whole copy is
complemented by a projection whose norm is majorized by $1+\delta'$
- a result which has recently been generalized by
N. Ozawa \cite{Oza-isom-embedd} to the category of operator spaces.

Secondly, can the $(m_n)$, $(a_n)$ and $(b_n)$ in Proposition
\ref{theo_A1} be arranged such that
$\norm{\phi_{m_n}-b_n \phi_{m_n} a_n}\gen0$ as $n\gen\infty$?
[Let us sketch in passing why this would generalize
Lemma \ref{lemma_A1}. If $(\phi_m)\subset A'$ is normalized and spans $\leins$
$r$-isomorphically then by James' distortion theorem there are blocks
$\psi_n=\sum_{i\in F_n}\lambda_i \phi_i$ spanning $\leins$
almost isometrically such that
$\sum_{i\in F_n}\betr{\lambda_i}\leq 1/r$.
Now, if there are appropriate $a_n, b_n\in A$
such that $\norm{\psi_n - b_n \psi_n a_n}\gen0$ (after passing,
if necessary, to an appropriate subsequence of $(\psi_n)$), one deduces that 
$\betr{\psi_n}(a_n)>\sqrt{1-\eps_n}$ with $0<\eps_n\gen0$.
Thus for each $n$ there is $i_n\in F_n$ such that
$\betr{\phi_{i_n}}(a_n) >(1-\eps_n)r^2$
because otherwise by \cite[III.4.7]{Tak}
one would have the contradiction
\bglst
\sqrt{1-\eps_n} &<& \betr{\psi_n}(a_n)
=\betr{\sum_{i\in F_n}\lambda_i \phi_i}(a_n)
\leq \Bigl(\sum_{F_n}\betr{\lambda_i}\,\norm{\phi_i}\Bigr)^{1/2}
    \Bigl(\sum_{F_n}\betr{\lambda_i}\,(\betr{\phi_i}(a_n^2))\Bigr)^{1/2}\\
&\leq& \frac{1}{r}\Bigl(\max_{F_n}\betr{\phi_i}(a_n^2)\Bigr)^{1/2}
\leq\sqrt{1-\eps_n}.\eglst
Similarly one obtains $\betr{\phi_{i_n}^*}(b_n) >(1-\eps_n)r^2$.]
\bigskip\\
Proposition \ref{theo_A1} follows immediately from Lemma \ref{lemma_A1}
(and \Ref{glA3_3} of Lemma \ref{lem_A3}) with $s=\eins=t$.
The technical part concerning $s,t$ is added because it might be usefull
for answering the second question just mentionned above.
\begin{lemma}\label{lemma_A1}
Let $A$ be a \Cstalg\ (unital or not), $r>0$, let $(\phi_m)$ be
a normalized sequence in $A'$ spanning $\leins$ $r$-isomorphically
that is such that
\bgl
r\sum\betr{\alpha_m}\leq\norm{\sum\alpha_m \phi_m }\leq
            \sum\betr{\alpha_m}
\falle (\alpha_{m})\aus {\C}.   \label{gl_leins_r1}
\egl
Then, given $\eps>0$, there are a sequence $(m_n)$ in $\N$ and
a sequence $(a_n)$ of pairwise orthogonal positive normalized
elements in $A$ and another sequence $(b_n)$ of pairwise
orthogonal positive normalized elements in $A$ such that
\bgl
\betr{\phi_{m_n}}(a_n) &>&
    (1-\eps)r^2
                            \label{gl_leins_r2}\\
\betr{\phi^{*}_{m_n}}(b_n) &>&
    (1-\eps)r^2
                            \label{gl_leins_r3}
\egl
for each $n\in\N$.\\
Moreover, if $s$ and $t$ are open projections in $A''$ such that $s$
($t$) majorizes the right (left) supports of all $\phi_m$
(that is $t\phi_m s=\phi_m$ for all $m\in\N$)
then one can obtain in addition that $a_n\in sA''s$ and
$b_n\in tA''t$.\\
Moreover, if the $\phi_m$ are selfadjoint one can obtain in addition
$\betr{\phi_{m_n}}(a_n)=\betr{\phi_{m_n}^*}(a_n)>(1-\eps)r$
instead of $\Ref{gl_leins_r2}$ and $\Ref{gl_leins_r3}$.
\end{lemma}
{\em Proof:}
First we suppose that $A$ is unital.\\
It is enough to construct a sequence $(p_n)$
of orthogonal open projections in $sA''s$ such that
\bgl
\betr{\phi_{m_n}}(p_n) &>&
    (1-\eps)r^2          \label{gl_leins_r2a}
\egl
for an appropriate subsequence $(\phi_{m_n})$
because then, by the definition of open projections,
for all $n\in\N$ positive elements $a_n\leq p_n$
can be choosen so to be pairwise orthogonal  (since the $p_n$ are)
and so to satisfy \Ref{gl_leins_r2}; finally, since
\Ref{gl_leins_r1} remains valid if $\phi_n^*$ is substituted for
$\phi_n$ the same reasoning that leads to \Ref{gl_leins_r2}
shows the existence of
a sequence $(b_n)$ in $tA''t$ as desired in \Ref{gl_leins_r3}.

Let $0<\eps<1$ and choose a sequence
$(\eps_{n})$ of positive numbers such that $\sum\eps_{n}=\eps$
and $\eps_{n}\leq \frac{3}{4}$ for all $n\in \N$.

By induction over $n=1,2,\ldots$ we construct a sequence
$(p_{n})$ of open projections in $sA''s$,
a sequence of indices
$(m_{n})$, a decreasing sequence
$(N_{n})$ of infinite subsets of $\N$, i.e.\ 
$\cdots \aus
N_{n+1}\aus N_{n}\aus \cdots \aus N_{1}\aus N_{0}=\N$,
such that we have for all $n\in\N$:
\bgl
\overline{p_{n}}&\in& sA''s
\label{gl74antibis}\\
\overline{p_{i}}\;\,\overline{p_{n}}&=&0
\falle i<n                                \label{gl74}\\
\overline{p_n}(\betr{\phi_{m}})   &<&
     \frac{1}{72}r^{2}\eps_{n}^{4} \;\;\;\;   \falle m\in N_{n}
\label{gl75}\\
p_n(\betr{\phi_{m_n}})    &>&
                                 r^2(1-\sum_{1}^{n}\eps_{i})
 \label{gl76}\\
m_{n}     &\in &       N_{n-1}.\label{gl77}
\egl
We start the induction with $n=1$.\\
Choose $j_{1}\in\N$ with $1/j_{1}<r^{2}\eps_{1}^{4}/72$.
For $j_{1}$ and $\eps_{1}/4$ Lemma \ref{lem_leins_endl}
yields a number
$\delta_{1}=\delta_{1}(j_{1},\eps_{1}/4)>0$ and without loss of
generality we assume
$\delta_{1}\leq \eps_{1}/4$.
By James' distortion theorem applied to \Ref{gl_leins_r1}
there are pairwise disjoint finite sets
$F_{k}^{(1)}\aus N_{0}=\N$,
a finite sequence $(\lambda_{i}^{(1)})_{i\in F_{k}^{(1)}}\aus\C$
and functionals
$\tau_{k}^{(1)}=\sum_{i\in F_{k}^{(1)}}\lambda_{i}^{(1)}\phi_{i}$ for
$k\in\N$, such that
\bgl
\sum_{F_{k}^{(1)}}\betr{\lambda_{i}^{(1)}}
&\leq& \frac{1}{r},
\label{gl71}\\
(1-\delta_{1})\sum_{k\geq 1}\betr{\alpha_{k}}
&\leq&
\norm{\sum_{k\geq 1}\alpha_{k}\tau_{k}^{(1)}}
\leq \sum_{k\geq 1}\betr{\alpha_{k}}
\falle (\alpha_{k})_{k\in N_{0}}\aus {\C}.
\label{gl72}
\egl
Again by Lemma \ref{lem_leins_endl}
there are pairwise orthogonal
open projections
$p_{k}^{(1)}\in sA''s$, $k\leq j_{1}$,
such that
\bgl
p_{k}^{(1)}(\betr{\tau_{k}^{(1)}})
  >(1-\eps_1/4)\norm{\tau_{k}^{(1)}}\falle k\leq j_{1},
\label{gl73}
\egl
and since
the projections can be chosen to have
or\-tho\-go\-nal closures in $sA''s$
we have
\bglst
\Bigl(\sum_{1}^{j_{1}}\overline{p_{k}^{(1)}}\Bigr)
(\betr{\phi_{m}})
\leq 1
\falle m\in N_{0}.
\eglst
Therefore there exist a $k_{1}\leq j_{1}$ 
and an infinite set $N_{1}\aus N_{0}$ such that
\bgl
\overline{p_{k_{1}}^{(1)}}    (\betr{\phi_{m}})
\leq \frac{1}{j_1}
< \frac{r^{2}\eps_{1}^{4}}{72} \falle m\in N_{1}.
\label{gl7a1}
\egl
Set $p_{1}=p_{k_{1}}^{(1)}$,
$\tau_{1}=\tau_{k_{1}}^{(1)}$, $F_{1}=F_{k_{1}}^{(1)}$.
Then \Ref{gl74antibis} holds for $n=1$.
Now we infer that
\bglst
p_{1}(\betr{\tau_{1}})
\stackrel{\Ref{gl73}}{>}
\norm{\tau_{1}}(1-\frac{\eps_{1}}{4})
\stackrel{\Ref{gl72}}{\geq}
 (1-\delta_{1})(1-\frac{\eps_{1}}{4})\geq
(1-\frac{\eps_{1}}{4})^2>\sqrt{1-\eps_1},
\eglst
which in turn yields the existence of an index
$m_{1}\in F_{1}\aus N_{0}$
as desired in \Ref{gl76} and \Ref{gl77}
for $n=1$, because otherwise we would
have
\bglst
p_{1}(\betr{\tau_{1}})
&=&
p_1\Bigl(\betr{\sum_{F_{1}}\lambda_{i}^{(1)}\phi_{i} }\Bigr)\\
& \stackrel{(*)}{\leq}&
\Bigl(\sum_{F_{1}}\norm{\lambda_i^{(1)}\phi_i}\Bigr)^{1/2}
\Bigl(\sum_{F_{1}}p_{1}
\bigl(\betr{\lambda_i^{(1)}\phi_{i}}\bigr)\Bigr)^{1/2}\\
&\leq&
\Bigl(\sum_{F_{1}}\betr{\lambda_i^{(1)}}\Bigr)^{1/2}
\Bigl(\max_{i\in F_{1}}p_{1}(\betr{\phi_{i}})
                  \sum_{F_{1}}\betr{\lambda_{i}^{(1)}}\Bigr)^{1/2}\\
&\stackrel{\Ref{gl71}}{\leq}&
\frac{1}{r}\bigl(r^2(1-\eps_{1})\bigr)^{1/2}=\sqrt{1-\eps_1}.
\eglst
Here inequality $(*)$ follows from \cite[III.4.7]{Tak}.
For (\ref{gl74}, n=1) nothing needs to be proved.
Inequality (\ref{gl75}, n=1) corresponds to
\Ref{gl7a1}. The first induction step is done.
\medskip\\
Induction step $n\gen n+1$:\\
Suppose $p_{k}$, $N_{k}$, $m_{k}$ to be constructed for $k\leq n$
according to         
\Ref{gl74} -- \Ref{gl77}.\\
Since the $\overline{p_{k}}$ are
orthogonal in $sA''s$, $\sum_{1}^{n}\overline{p_{k}}$
is closed by \cite[Th.\ II.7]{Ake69}.
Therefore $s_{n}=s-\sum_{1}^{n}\overline{p_{k}}\in sA''s$ is open.
Set $\tilde{\phi}_{m}=\phi_m s_n$.\\
{\sc Claim}:\\
The normalized functionals
$\bigl(\frac{\tilde{\phi}_{m}}{\norm{\tilde{\phi}_{m}}}
             \bigr)_{m\in N_{n}}$
form an $\leins$-basis with
\bgl
r\Bigl(1-\sum_{1}^{n}\eps_{i}^{2}\Bigr)^{1/2}
\sum_{m\in N_{n}}\betr{\alpha_{m}}
& \leq  &
\norm{\sum_{m\in N_{n}}\alpha_{m}
\frac{\tilde{\phi}_{m}}{\norm{\tilde{\phi}_{m}}}}\nonumber\\
&\leq &
\sum_{m\in N_{n}}\betr{\alpha_{m}}
\falle (\alpha_{m})_{m\in N_{n}}\aus {\C}.
\label{gl730}
\egl
Set $\eta=\frac{r^2}{72}\sum_1^{n}\eps_{i}^{4}$.
Then
\bgl
(s-s_{n})(\betr{\phi_{m}})=
\Bigl(\sum_{1}^{n}
\overline{p_{k}}\Bigr)(\betr{\phi_{m}})
\stackrel{\Ref{gl75}}{<}                         \eta
\falle m\in N_{n},      \label{gl_leins10a}
\egl
thus since $s(\betr{\phi_m}) = \norm{\phi_m}=1$
\bgl
\norm{\phi_m s_n -\phi_m}
&\stackrel{\Ref{glA3_3}}{\leq}&
\betr{2(\norm{\phi_m}-s_n(\betr{\phi_m}))}^{1/2}\nonumber\\
&=&
\bigl(2\sum_1^n \overline{p}_k (\betr{\phi_m})\bigr)^{1/2}
\stackrel{\Ref{gl_leins10a}}{<}
\sqrt{2\eta}
\falle m\in N_{n}; \label{gl735}
\egl
further we note that for all $m\in N_{n}$
\bgl
\norm{\tilde{\phi}_{m}}&\leq& \norm{\phi_{m}}=1,
\label{gl7b1}\\
0\leq 1-\norm{\tilde{\phi}_{m}}&=&
\norm{\phi_{m}}-\norm{\tilde{\phi}_{m}}
\leq \norm{\phi_{m}-\tilde{\phi}_{m}}
\stackrel{\Ref{gl735}}{\leq}
\sqrt{2\eta}
\label{gl7b2}\\
{\norm{\tilde{\phi}_{m}}}
&\stackrel{\Ref{gl7b2}}{\geq}&
1-\sqrt{2\eta},
\label{gl7b3}
\egl
hence
\bgl
\norm{\frac{\tilde{\phi}_{m}}{\norm{\tilde{\phi}_{m}}}-\tilde{\phi}_{m}}
&\stackrel{\Ref{gl7b1}}{\leq}&
\frac{1}{\norm{\tilde{\phi}_{m}}}-1 =
\Bigl(1-\norm{\tilde{\phi}_{m}}\Bigr)\frac{1}{\norm{\tilde{\phi}_{m}}}
\nonumber\\
&\stackrel{\Ref{gl7b2}\,\Ref{gl7b3}}{\leq}&
\frac{\sqrt{2\eta}}{1-\sqrt{2\eta}}
\label{gl7b4}
\egl
and
\bgl
\norm{\frac{\tilde{\phi}_{m}}{\norm{\tilde{\phi}_{m}}}-\phi_{m}}
&\leq&
\norm{\frac{\tilde{\phi}_{m}}{\norm{\tilde{\phi}_{m}}}-\tilde{\phi}_{m}}+
\norm{\tilde{\phi}_{m}-\phi_{m}}\nonumber\\
&\stackrel{\Ref{gl7b4}\,\Ref{gl735}}{\leq}&
\sqrt{2\eta}\Bigl(1+\frac{1}{1-\sqrt{2\eta}}\Bigr)
<3\sqrt{2\eta} \label{gl7351}
\egl
because $\eps<1$, $r\leq1$, thus $\sqrt{2\eta}<1/2$.
Then \Ref{gl730} follows from
\bglst
\norm{\sum_{m\in N_{n}}\alpha_{m}
\frac{\tilde{\phi}_{m}}{\norm{\tilde{\phi}_{m}}}}
&\geq&
\norm{\sum_{N_{n}}\alpha_{m}\phi_{m}} - \norm{\sum_{N_{n}}\alpha_{m}
\frac{\tilde{\phi}_{m}}{\norm{\tilde{\phi}_{m}}}-
\sum_{N_{n}}\alpha_{m}\phi_{m}}\\
&\stackrel{\Ref{gl_leins_r1}}{\geq}&
r\Bigl(1-\sup_{N_{n}}
\norm{\frac{\tilde{\phi}_{m}}{\norm{\tilde{\phi}_{m}}}
-\phi_{m}}\Bigr)\sum_{N_{n}}\betr{\alpha_{m}}\\
&\stackrel{\Ref{gl7351}}{\geq}&
r(1-3\sqrt{2\eta})\sum_{N_{n}}\betr{\alpha_{m}}\\
&=&
r\Bigl(1-\frac{1}{2}(\sum_1^n \eps_i^4)^{1/2}\Bigr)
\sum_{N_{n}}\betr{\alpha_{m}}\\
&>&r\Bigl(1-(\sum_1^n \eps_i^4)^{1/2}\Bigr)^{1/2}
    \sum_{N_{n}}\betr{\alpha_{m}}\\
&>&
r\Bigl(1-\sum_{1}^{n}\eps_{i}^{2}\Bigr)^{1/2}
\sum_{m\in N_{n}}\betr{\alpha_{m}}
\eglst
and the {\sc Claim} is established.

Choose a number $j_{n+1}\in\N$ such that
$1/j_{n+1}<r^{2}\eps_{n+1}^{2}/4$.
Further choose a number
$\delta_{n+1}=\delta_{n+1}(j_{n+1},\eps_{n+1}^{2}/4)>0$
according to Lemma \ref{lem_leins_endl} and such that moreover
$\delta_{n+1}\leq \eps_{n+1}^{2}/4$.
Now we apply  James' distortion theorem.
By \Ref{gl730} there are
pairwise disjoint finite sets
$F_{k}^{(n+1)}\aus N_{n}$,
a finite sequence $(\lambda_{i}^{(n+1)})_{i\in F_{k}^{(n+1)}}\aus\C$
and functionals
$\tau_{k}^{(n+1)}=\sum_{i\in F_{k}^{(n+1)}}\lambda_{i}^{(n+1)}
\frac{\tilde{\phi}_{i}}{\norm{\tilde{\phi}_{i}}}$ for each $k\in\N$ 
such that
\bgl
\sum_{i\in F_{k}^{(n+1)}}\betr{\lambda_{i}^{(n+1)}}
&\leq&
 \frac{1}{r(1-\sum_{1}^{n}\eps_{i}^{2})^{1/2}} \falle k\in {\N},
\label{gl71n}\\
(1-\delta_{n+1})\sum_{k\geq 1}\betr{\alpha_{k}}
&\leq&
\norm{\sum_{k\geq 1}\alpha_{k}\tau_{k}^{(n+1)}s_n}
\leq \sum_{k\geq 1}\betr{\alpha_{k}}.
\label{gl72n}
\egl
Again by Lemma \ref{lem_leins_endl},
applied to the open projections $s_n$ and $\eins$,
to the functionals
$\tau_{k}^{(n+1)}\in A'$, and to \Ref{gl72n},
there exist open projections
$p_{k}^{(n+1)}\in s_n A'' s_n$, $k\leq j_{n+1}$,
with pairwise orthogonal closures in $s_n A''s_n$ such that
\bgl
p_{k}^{(n+1)}(\betr{\tau_{k}^{(n+1)}})>
\norm{\tau_{k}^{(n+1)}} \Bigl(1-\frac{\eps_{n+1}^{2}}{4}\Bigr)
\label{gl73n}
\egl
for $k\leq j_{n+1}$.
Since the projections have orthogonal closures
we have
\bglst
\Bigl(\sum_{1}^{j_{n+1}}
\overline{p_{k}^{(n+1)}}
\Bigr)(\betr{\phi_{m}})
\leq 1
\;\,\,\,\falle m\in N_{n}.
\eglst
Therefore there exist
an index $k_{n+1}\leq j_{n+1}$ and an infinite subset
$N_{n+1}\aus N_{n}$ such that
\bglst
\overline{p_{k_{n+1}}^{(n+1)}}(\betr{\phi_{m}})
\leq
      \frac{1}{j_{n+1}}
<
   \frac{r^{2}\eps_{n+1}^{4}}{72} \falle m\in N_{n+1}.
\eglst
Set $p_{n+1}=p_{k_{n+1}}^{(n+1)}$,
$\tau_{n+1}=\tau_{k_{n+1}}^{(n+1)}$, $F_{n+1}=F_{k_{n+1}}^{(n+1)}$.
Then (\ref{gl74antibis}), \Ref{gl74} and \Ref{gl75} hold for $n+1$.
Now we infer that
\bgl
p_{n+1}(\betr{\tau_{n+1}})
&\stackrel{\Ref{gl73n}}{>}&
\norm{\tau_{n+1}}
               \Bigl(1-\frac{\eps_{n+1}^{2}}{4}\Bigr) \nonumber\\
&\stackrel{\Ref{gl72n}}{\geq}&
 (1-\delta_{n+1})(1-\frac{\eps_{n+1}^{2}}{4})
    >
\sqrt{1-\eps_{n+1}^{2}},\label{gl74n}
\egl
hence there is an index $m_{n+1}\in F_{n+1}\aus N_{n}$
as desired in \Ref{gl77}
such that
\bgl
p_{n+1}\bigl(\frac{\tilde{\phi}_{m_{n+1}}}
{\norm{\tilde{\phi}_{{m_{n+1}}}}}\bigr)>
r^2(1-\eps_{n+1}^{2})\Bigl(1-\sum_{1}^{n}\eps_{i}^{2}\Bigr),
\label{gl7b8}
\egl
because otherwise the following estimates would contradict
\Ref{gl74n}:
\bglst
p_{n+1}(\betr{\tau_{n+1}})
&=&
p_{n+1}\bigl(\betr{\sum_{F_{n+1}}\lambda_{i}^{(n+1)}
    \frac{\tilde{\phi_i}}{\norm{\tilde{\phi_i}}}}\bigr)\\
&\stackrel{(*)}{\leq}&
\bigl(\sum_{F_{n+1}}\norm{\lambda_i^{(n+1)}\frac{\tilde{\phi_i}}{\norm{\tilde{\phi_i}}}}\bigr)^{1/2}
\bigl(\sum_{F_{n+1}}p_{n+1}(\betr{\lambda_i^{(n+1)}
\frac{\tilde{\phi_i}}{\norm{\tilde{\phi_i}}}})\bigr)^{1/2}                         \\
&\leq&
\bigl(\sum_{F_{n+1}}\betr{\lambda_i^{(n+1)}}\bigr)^{1/2}
\bigl(\max_{i\in F_{n+1}}p_{n+1}(\betr{\frac{\tilde{\phi_i}}{\norm{\tilde{\phi_i}}}})
                  \sum_{F_{n+1}}\betr{\lambda_{i}^{(n+1)}}\bigr)^{1/2}\\
&\stackrel{\Ref{gl71n}}{\leq}&
\frac{1}{r(1-\sum_1^n \eps_i^2)^{1/2}}
\bigl(r^2(1-\eps_{n+1}^2)(1-\sum_1^n \eps_i^2)\bigr)^{1/2}
=
\sqrt{1-\eps_{n+1}^2}.
\eglst
Here inequality $(*)$ follows from \cite[III.4.7]{Tak}.
Note that for a functional $\phi$ with polar decomposition 
$\phi=u\betr{\phi}$ one has
$\betr{\phi s}=\betr{u\,\betr{\phi}\, s} \leq\betr{\phi}$
by \cite[III.4.9]{Tak} which explains inequality $(**)$ below;
now (\ref{gl76}, $n+1$) follows from
\bglst
p_{n+1}(\betr{\phi_{m_{n+1}}})
&\stackrel{(**)}{\geq}&
 p_{n+1}(\betr{\phi_{m_{n+1}}s_n})
\\
&\stackrel{\Ref{gl7b8}}{>}&
\norm{\phi_{m_{n+1}}s_n}r^2
(1-\eps_{n+1}^{2})\Bigl(1-\sum_{1}^{n}\eps_{i}^{2}\Bigr)\\
&\stackrel{\Ref{gl7b3}}{\geq}&
\Bigl(1-\sqrt{2\eta} \Bigr)
r^2(1-\eps_{n+1}^{2}) \Bigl(1-\sum_{1}^{n}\eps_{i}^{2}\Bigr)\\
&>&
r^2 \Bigl(1-\sum_{1}^{n+1}\eps_{i}\Bigr)
\eglst
where the last inequality follows from the following
completely elementary
estimates:
\bglst
\lefteqn{
\Bigl(1-\sqrt{2\eta} \Bigr)
(1-\eps_{n+1}^{2})\Bigl(1-\sum_{1}^{n}\eps_{i}^{2}\Bigr)=
         }
\\&=&
\Bigl(1-\frac{r}{6}(\sum_{1}^{n}\eps_{i}^{4})^{1/2}\Bigr)
(1-\eps_{n+1}^{2})\Bigl(1-\sum_{1}^{n}\eps_{i}^{2}\Bigr)\\
&>&
\Bigl(1-\frac{r}{3}\sum_{1}^{n}\eps_{i}^{2}\Bigr)
(1-\eps_{n+1}^{2})\Bigl(1-\sum_{1}^{n}\eps_{i}^{2}\Bigr)\\
&=&
\Bigl(1-\sum_{1}^{n+1}\eps_{i}^{2}\Bigr)-
   \frac{r}{3}\Bigl[1-(1-\eps_{n+1}^{2})\sum_{1}^{n}\eps_{i}^{2}
-\bigl(1+\frac{3}{r}\bigr)\eps_{n+1}^{2} \Bigr]
\sum_{1}^{n}\eps_{i}^{2}\\
&>&
\Bigl(1-\sum_{1}^{n+1}\eps_{i}^{2}\Bigr)
-\frac{r}{3}\sum_{1}^{n}\eps_{i}^{2}
>
\Bigl(1-\sum_{1}^{n+1}\eps_{i}^{2}\Bigr)
-\frac{1}{3}\sum_{1}^{n+1}\eps_{i}^{2}\\
&=&
\Bigl(1-\sum_{1}^{n+1}\eps_{i}\Bigr)+
\sum_{1}^{n+1}\eps_{i}(1-\eps_{i}-\frac{\eps_{i}}{3})
\geq
1-\sum_{1}^{n+1}\eps_{i}
\eglst
since we assumed $\eps_{i}\leq \frac{3}{4}$ for all $i\in\N$.
Thus (\ref{gl76}, $n+1$) is proved.
This ends the induction and the proof if $A$ is unital.\\
If $A$ is not unital we consider its unitisation $A_{\eins}=A\oplus \C \eins$
on which multiplication is defined by $(a,\lambda)(b,\mu) =
(ab+\lambda b+\mu a, \lambda\mu)$. Note that if $\hat{\phi}\in A_{\eins}'$
is a norm preserving Hahn-Banach extension of $\phi\in A'$ then 
$\betr{\hat{\phi}}\in A_{\eins}'$ is an extension of $\betr{\phi}\in A'$.
[To see this let $u_n\in A$, $\norm{u_n}\leq 1$ be such that 
$\norm{\hat{\phi}}=\norm{\phi}=\lim \phi(u_n)$. Then both 
$\,\norm{\betr{\phi} - u_n \phi}\gen0$ and
$\,\norm{\betr{\hat{\phi}} - (u_n,0) \hat{\phi}}\gen0$ by
\Ref{glA3_2} of Lemma \ref{lem_A3}.
Hence, $\betr{\hat{\phi}}((x,0))=\lim \hat{\phi}((x,0)(u_n,0))
=\lim \hat{\phi}((xu_n,0)) =\lim \phi(xu_n)
=\betr{\phi}(x)$ for any $x=(x,0)\in A$.]
Let $\hat{\phi}_m$ be a norm preserving Hahn-Banach extension of $\phi_m$.
Then \Ref{gl_leins_r1} remains valid and by what has been proved
for the unital case we get normalized pairwise orthogonal positive
$(a_n,\lambda_n)\in A_{\eins}$ such that $(a_n,\lambda_n)\leq s$ and
$\betr{\hat{\phi}_{m_n}}((a_n,\lambda_n)) > (1-\eps)r^2$ for an
appropriate subsequence $(\hat{\phi}_{m_n})$.
Since the $(a_n,\lambda_n)$ are pairwise orthogonal all
(but possibly one) of them have $\lambda_n=0$ whence 
$\betr{\phi_{m_n}}(a_n)=
  \betr{\hat{\phi}_{m_n}}((a_n, 0)) > (1-\eps)r^2$.
Likewise, we get \Ref{gl_leins_r3}.\\
The last statement concerning selfadjoint functionals
has been discussed in the beginning of this section.
\medskip
\ebew\bigskip\\
{\sc Acknowledgement}
I thank Dirk Werner for several helpful discussions.

\bigskip
Hermann Pfitzner\\
Universit\'e d'Orl\'eans\\
BP 6759\\
F-45067 Orl\'eans Cedex 2\\
France\\
e-mail: pfitzner@labomath.univ-orleans.fr
\end{document}